# ON FORMULAE DECOUPLING THE TOTAL VARIATION OF BV FUNCTIONS

AUGUSTO C. PONCE AND DANIEL SPECTOR

*To Nicola Fusco, a master of BV, on the occasion of his sixtieth birthday*

ABSTRACT. In this paper we prove several formulae that enable one to capture the singular portion of the measure derivative of a function of bounded variation as a limit of non-local functionals. One special case shows that rescalings of the fractional Laplacian of a function $u \in SBV$ converge strictly to the singular portion of $Du$.

## 1. INTRODUCTION AND MAIN RESULTS

Let $\Omega \subset \mathbb{R}^N$ be an open, bounded and smooth subset (or all of $\mathbb{R}^N$) and $u$ be a function of bounded variation in $\Omega$, i.e. $u \in L^1(\Omega)$ and its distributional derivative $Du$ is a Radon measure with finite total variation

$$(1.1) \quad |Du|(\Omega) := \sup_{\substack{\Phi \in C_c^1(\Omega; \mathbb{R}^N), \\ \|\Phi\|_{L^\infty(\Omega; \mathbb{R}^N)} \leq 1}} \int_\Omega u \operatorname{div} \Phi.$$

The space of such functions $BV(\Omega)$ contains the Sobolev space $W^{1,1}(\Omega)$, with strict inclusion, since in general for $u \in BV(\Omega)$ one has the decomposition of the measure

$$Du = \nabla u \, \mathcal{L}^N + D^s u,$$

where $\nabla u$ is the Radon-Nikodym density of $Du$ with respect to the $N$-dimensional Lebesgue measure $\mathcal{L}^N$ and $D^s u$ is singular with respect to $\mathcal{L}^N$. In particular, $u \in W^{1,1}(\Omega)$ precisely when $D^s u \equiv 0$.

While the computation of the total variation (1.1) through the theory of distributions is classical, in recent years there has been an interest in its — and other related energies — approximation through the asymptotics of non-local functionals [1, 4–7, 9, 11, 12, 16–18, 20, 22–26]. For example, Bourgain, Brezis and Mironescu [7] had shown that for $u \in W^{1,p}(\Omega)$, $1 \leq p < +\infty$, one has

$$(1.2) \quad \lim_{\alpha \to 1^-} (1-\alpha) \int_\Omega \int_\Omega \frac{|u(x) - u(y)|^p}{|x-y|^{N+\alpha p}} \, dy \, dx = C_{p,N} \int_\Omega |\nabla u|^p.$$

Here, for $p = 1$,

$$C_{1,N} = \int_{\mathbb{S}^{N-1}} |e \cdot v| \, d\mathcal{H}^{N-1}(v) = 2\omega_{N-1},$$

where $e \in \mathbb{R}^N$ is any unit vector and $\omega_{N-1}$ is the volume of the unit ball in $\mathbb{R}^{N-1}$. Formula (1.2) expresses the fact that appropriately scaled Gagliardo semi-norms tend to the total variation as the differential parameter tends to one. Their result includes the case $u \in W^{1,1}(\Omega)$, while for general $u \in BV(\Omega)$ it was proved by Dávila [14, Theorem 1] that

$$(1.3) \quad \lim_{\alpha \to 1^-} (1-\alpha) \int_\Omega \int_\Omega \frac{|u(x) - u(y)|}{|x-y|^{N+\alpha}} \, dy \, dx = C_{1,N}|Du|(\Omega).$$

The Gagliardo semi-norms can be thought heuristically as the energy of derivatives of fractional order $\alpha \in (0,1)$, while the space of $p$-integrable functions for which they are finite coincides with the $W^{\alpha,p}$ space obtained by *real interpolation* of $L^p$ and $W^{1,p}$ with parameter $\alpha$. Indeed, it was subsequently shown by Milman [24] that the convergence (1.2) can be alternatively deduced from the theory of interpolation.

In the *complex interpolation* between $L^p$ and $W^{1,p}$ when $\Omega = \mathbb{R}^N$, one classically sees the fractional Laplacian

$$(-\Delta)^{\alpha/2} u(x) := c_\alpha \int_{\mathbb{R}^N} \frac{u(x) - u(y)}{|x-y|^{N+\alpha}} \, dy$$

filling the role of the differential object of order $\alpha \in (0,1)$, and even in a more extended range, though with a different formula. The constant

$$c_\alpha = \frac{2^{\alpha-1} \alpha \Gamma(\frac{N+\alpha}{2})}{\pi^{\frac{N}{2}} \Gamma(\frac{2-\alpha}{2})}$$

ensures that the Fourier transform of the fractional Laplacian satisfy $\widehat{(-\Delta)^{\alpha/2} u}(\xi) = (2\pi|\xi|)^\alpha \hat{u}(\xi)$.

As integer powers of the Laplacian can also be understood from the framework of the functional calculus, the natural limits here are $\alpha = 0$ and $\alpha = 2$. This observation is supported by the trivial convergence, as $\alpha$ tends to 1, of the energies

$$(1.4) \qquad \lim_{\alpha \to 1^-} (1-\alpha) \int_{\mathbb{R}^N} \left| \int_{\mathbb{R}^N} \frac{u(x) - u(y)}{|x-y|^{N+\alpha}} \, dy \right|^p \, dx = 0,$$

whenever $u \in W^{1,p}(\mathbb{R}^N)$ and $p \geq 1$. The case $p > 1$ follows from the $L^p$ boundedness of the Riesz transforms that yields the equivalence between the quantities $\|\nabla u\|_{L^p(\mathbb{R}^N)}$ and $\|(-\Delta)^{1/2} u\|_{L^p(\mathbb{R}^N)}$; see e.g. [19, Lemma 3.6]. The case $p = 1$ requires a different argument that can be found in [27, 28], based on the approximation of $u$ by smooth functions. Alternatively, one can still obtain the expected energy, in the spirit of (1.2), via a non-local gradient approach [23].

What is perhaps surprising is that the convergence (1.4) is no longer true for general $u \in BV(\Omega)$ and $p = 1$. For example, in the case where $u = \chi_A$ is the characteristic function of a set of finite perimeter, Dávila's formula (1.3) and the straightforward integral identity

$$(1.5) \qquad \int_{\mathbb{R}^N} \left| \int_{\mathbb{R}^N} \frac{\chi_A(x) - \chi_A(y)}{|x-y|^{N+\alpha}} \, dy \right| \, dx = \int_{\mathbb{R}^N} \int_{\mathbb{R}^N} \frac{|\chi_A(x) - \chi_A(y)|}{|x-y|^{N+\alpha}} \, dy \, dx,$$

imply that

$$(1.6) \qquad \lim_{\alpha \to 1^-} (1-\alpha) \int_{\mathbb{R}^N} \left| \int_{\mathbb{R}^N} \frac{\chi_A(x) - \chi_A(y)}{|x-y|^{N+\alpha}} \, dy \right| \, dx = C_{1,N} |D^s \chi_A|(\mathbb{R}^N).$$

Here we observe that the measure $D\chi_A$ is singular with respect to $\mathcal{L}^N$, and thus $D\chi_A = D^s \chi_A$. More generally, in this paper we are interested in several formulae for the asymptotics of functionals capturing the *singular portion* of the measure derivative of $u \in BV(\Omega)$. Let us first restrict our attention to the setting of special functions of bounded variation $SBV(\Omega)$. This space has been introduced by De Giorgi and Ambrosio [15], and consists of those $BV$ functions whose singular part of the derivative is supported by an $(N-1)$-dimensional rectifiable set; see Section 3 below. Then our first result is

**Theorem 1.1.** *For every $u \in SBV(\Omega)$, we have*

$$\lim_{\epsilon \to 0} \int_\Omega \left| \int_\Omega \frac{u(x) - u(y)}{|x-y|} \rho_\epsilon(x-y) \, dy \right| \, dx = K_{1,N} |D^s u|(\Omega),$$

where $K_{1,N} = C_{1,N}/\mathcal{H}^{N-1}(\mathbb{S}^{N-1})$ and $\mathcal{H}^{N-1}(\mathbb{S}^{N-1})$ denotes the $(N-1)$-dimensional Hausdorff measure of the unit sphere $\mathbb{S}^{N-1}$.

We assume throughout the paper that $(\rho_\epsilon)_{\epsilon>0} \subset L^1(\mathbb{R}^N)$ is a family of non-negative radial functions that satisfies

$$\text{(1.7)} \qquad \lim_{\epsilon \to 0} \int_{|h|<\delta} \rho_\epsilon(h) \, dh = 1,$$

$$\text{(1.8)} \qquad \lim_{\epsilon \to 0} \int_{|h|>\delta} \rho_\epsilon(h) \, dh = 0,$$

for every $\delta > 0$. Taking in particular

$$\rho_\epsilon(h) = \frac{\epsilon}{\mathcal{H}^{N-1}(\mathbb{S}^{N-1})} \frac{\chi_{B_R(0)}(h)}{|h|^{N-\epsilon}},$$

for any fixed $R > 0$, one finds an analogue to the convergence (1.2) for $p = 1$. Moreover, as $c_1 C_{1,N} = 2/\pi$, one can show by a slight variation of the argument of Theorem 1.1 the following

**Corollary 1.2.** *For every $u \in SBV(\mathbb{R}^N)$, we have*

$$\lim_{\alpha \to 1^-} (1-\alpha) \|(-\Delta)^{\alpha/2} u\|_{L^1(\mathbb{R}^N)} = \frac{2}{\pi} |D^s u|(\mathbb{R}^N).$$

In Section 4 below, we show that the standard Cantor function admits a positive lower bound for some families $(\rho_\epsilon)_{\epsilon>0}$, including $\rho_\epsilon = \chi_{(-\epsilon,\epsilon)}/2\epsilon$. This motivates

**Open Problem 1.3.** *Do these identities hold for every $u \in BV(\mathbb{R}^N)$?*

The answer is not clear. For example, Fusco, Moscariello and Sbordone [18] have recently studied a related functional introduced by Ambrosio, Bourgain, Brezis and Figalli [1, Section 4.2]. In their case, the contribution of the singular part of the derivative cannot be expressed solely in terms of the total mass $|D^s u|(\mathbb{R}^N)$ [18, Example 2.2 and Theorem 3.3]. In Section 5, we discuss some connections between these various functionals.

We next turn our attention to the entire space of functions of bounded variation. In this broader setting, we have the following result, an extension of the recent Taylor characterization by the second author [27, 28] to this regime, which bears some analogy with formula (1.6):

**Theorem 1.4.** *For every $u \in BV(\Omega)$, we have that*

$$\lim_{\epsilon \to 0} \int_\Omega \int_\Omega \frac{|u(x) - u(y) - \nabla u(x) \cdot (x-y)|}{|x-y|} \rho_\epsilon(x-y) \, dy \, dx = K_{1,N} |D^s u|(\Omega).$$

**Remark 1.5.** *Theorem 1.4 combined with the results of [7] or [23] yields the following characterization of the Sobolev space $W^{1,1}$: $u \in W^{1,1}(\Omega)$ if and only if*

$$\lim_{\epsilon \to 0} \int_\Omega \int_\Omega \frac{|u(x) - u(y) - F(x) \cdot (x-y)|}{|x-y|} \rho_\epsilon(x-y) \, dy \, dx = 0,$$

*for some function $F \in L^1(\Omega; \mathbb{R}^N)$, and in that case $F = \nabla u$ almost everywhere in $\Omega$.*

In fact we prove slightly stronger results than the convergence of the energies in Theorems 1.1 and 1.4 and Corollary 1.2, the so-called strict convergence of the measures defined by the integrands; see Propositions 2.1 and 3.2 below. The organization of the paper is as follows: as we utilize the result of Theorem 1.4 in the proof of Theorem 1.1, we first prove Theorem 1.4 in Section 2. In Section 3 we prove Theorem 1.1 and Corollary 1.2. Next in Section 4 we give the example of the non-trivial lower bound for the Cantor function. Finally, in Section 5 we discuss some of the connections between the functionals we consider and those of Bourgain, Brezis and Mironescu [7], Ambrosio, Bourgain, Brezis and Figalli [1] and Fusco, Moscariello and Sbordone [18].

## 2. Proof of Theorem 1.4

Given $u \in BV(\mathbb{R}^N)$, define the function $R_\epsilon^1 u : \mathbb{R}^N \to \mathbb{R}$ by

$$R_\epsilon^1 u(x) = \int_{\mathbb{R}^N} \frac{|u(x+h) - u(x) - \nabla u(x) \cdot h|}{|h|} \rho_\epsilon(h) \, \mathrm{d}h.$$

We are interested in the convergence of $\|R_\epsilon^1 u\|_{L^1(\mathbb{R}^N)}$ as $\epsilon$ tends to zero. Observe that if $u \in C_c^\infty(\mathbb{R}^N)$, then

$$R_\epsilon^1 u \to 0 \quad \text{in } L^1(\mathbb{R}^N).$$

This property still holds if $u \in C^\infty(\mathbb{R}^N)$ and $\nabla u \in L^1(\mathbb{R}^N; \mathbb{R}^N)$.

**Proposition 2.1.** *Let $u \in BV(\mathbb{R}^N)$. For every bounded continuous function $\varphi : \mathbb{R}^N \to \mathbb{R}$, we have*

$$\lim_{\epsilon \to 0} \int_{\mathbb{R}^N} R_\epsilon^1 u \, \varphi = K_{1,N} \int_{\mathbb{R}^N} \varphi \, \mathrm{d}|D^s u|.$$

We first prove a couple of lemmas. Given $u \in BV(\mathbb{R}^N)$ and a family of mollifiers $(\psi_\delta)_{\delta > 0}$ in $C_c^\infty(\mathbb{R}^N)$, we use the notation $u_\delta = u * \psi_\delta$, $(\nabla u)_\delta = \nabla u * \psi_\delta$, ...

**Lemma 2.2.** *For every $u \in BV(\mathbb{R}^N)$, we have*

$$K_{1,N} |D^s u|(\mathbb{R}^N) \leq \liminf_{\epsilon \to 0} \|R_\epsilon^1 u\|_{L^1(\mathbb{R}^N)}.$$

*Proof.* By comparison of integrals and Fubini's theorem, for every $x \in \mathbb{R}^N$ we have

$$\int_{\mathbb{R}^N} \frac{|u_\delta(x+h) - u_\delta(x) - (\nabla u)_\delta(x) \cdot h|}{|h|} \rho_\epsilon(h) \, \mathrm{d}h \leq (R_\epsilon^1 u * \psi_\delta)(x) = (R_\epsilon^1 u)_\delta(x).$$

Integration with respect to $x$ yields

$$\int_{\mathbb{R}^N} \int_{\mathbb{R}^N} \frac{|u_\delta(x+h) - u_\delta(x) - (\nabla u)_\delta(x) \cdot h|}{|h|} \rho_\epsilon(h) \, \mathrm{d}h \, \mathrm{d}x \leq \int_{\mathbb{R}^N} (R_\epsilon^1 u)_\delta(x) \, \mathrm{d}x = \|R_\epsilon^1 u\|_{L^1(\mathbb{R}^N)}.$$

On the other hand, since

$$K_{1,N} |(D^s u)_\delta(x)| = \int_{\mathbb{R}^N} \left| (D^s u)_\delta(x) \cdot \frac{h}{|h|} \right| \rho_\epsilon(h) \, \mathrm{d}h$$

and

$$(\nabla u)_\delta + (D^s u)_\delta = \nabla(u_\delta),$$

by the triangle inequality we have

$$\left| \int_{\mathbb{R}^N} \frac{|u_\delta(x+h) - u_\delta(x) - (\nabla u)_\delta(x) \cdot h|}{|h|} \rho_\epsilon(h) \, \mathrm{d}h - K_{1,N} |(D^s u)_\delta(x)| \right| \leq R_\epsilon^1(u_\delta)(x).$$

Hence,

$$K_{1,N} \int_{\mathbb{R}^N} |(D^s u)_\delta(x)| \, \mathrm{d}x \leq \int_{\mathbb{R}^N} R_\epsilon^1(u_\delta)(x) \, \mathrm{d}x + \int_{\mathbb{R}^N} \int_{\mathbb{R}^N} \frac{|u_\delta(x+h) - u_\delta(x) - (\nabla u)_\delta(x) \cdot h|}{|h|} \rho_\epsilon(h) \, \mathrm{d}h \, \mathrm{d}x$$
$$\leq \|R_\epsilon^1(u_\delta)\|_{L^1(\mathbb{R}^N)} + \|R_\epsilon^1 u\|_{L^1(\mathbb{R}^N)}.$$

Since $u_\delta \in C^\infty(\mathbb{R}^N)$ and $\nabla(u_\delta) \in L^1(\mathbb{R}^N; \mathbb{R}^N)$ for every $\delta > 0$, the first term in the right-hand side converges to zero as $\epsilon$ tends to zero, and we get

$$K_{1,N} \int_{\mathbb{R}^N} |(D^s u)_\delta(x)| \, \mathrm{d}x \leq \liminf_{\epsilon \to 0} \|R_\epsilon^1 u\|_{L^1(\mathbb{R}^N)}.$$

Letting $\delta$ tend to zero, the conclusion follows from the lower semicontinuity of the norm. □

**Lemma 2.3.** *Let $u \in BV(\mathbb{R}^N)$. For every non-negative bounded continuous function $\varphi : \mathbb{R}^N \to \mathbb{R}$, we have*

$$\limsup_{\epsilon \to 0} \int_{\mathbb{R}^N} R_\epsilon^1 u \, \varphi \leq K_{1,N} \int_{\mathbb{R}^N} \varphi \, d|D^s u|.$$

*Proof.* We prove that

$$\begin{aligned}(2.1) \quad \int_{\mathbb{R}^N} R_\epsilon^1 u \, \varphi &\leq K_{1,N} \int_{\mathbb{R}^N} \varphi \, d|D^s u| \\ &+ \int_0^1 \int_{\mathbb{R}^N} \left( \int_{\mathbb{R}^N} |\varphi(x+th) - \varphi(x)| \rho_\epsilon(h) \, dh \right) d|D^s u|(x) \, dt \\ &+ \|\varphi\|_{L^\infty(\mathbb{R}^N)} \int_0^1 \int_{\mathbb{R}^N} \left( \int_{\mathbb{R}^N} |\nabla u(x+th) - \nabla u(x)| \, dx \right) \rho_\epsilon(h) \, dh \, dt.\end{aligned}$$

Since $\varphi$ is bounded and continuous and $\nabla u \in L^1(\mathbb{R}^N; \mathbb{R}^N)$, the second and third terms converge to zero as $\epsilon$ tends to zero. It thus suffices to establish this estimate. For this purpose, we proceed by approximation using the functions $u_\delta$. By the Fundamental Theorem of Calculus and the triangle inequality, we have that

$$|u_\delta(x+h) - u_\delta(x) - (\nabla u)_\delta(x) \cdot h| \leq \int_0^1 |(D^s u)_\delta(x+th) \cdot h| \, dt + \int_0^1 |(\nabla u)_\delta(x+th) - (\nabla u)_\delta(x)| \cdot |h| \, dt.$$

By the triangle inequality and the change of variables $z = x + th$, we also have that

$$\begin{aligned}&\int_{\mathbb{R}^N} |(D^s u)_\delta(x+th) \cdot h| \varphi(x) \, dx \\ &\leq \int_{\mathbb{R}^N} |(D^s u)_\delta(x+th) \cdot h| \varphi(x+th) \, dx + \int_{\mathbb{R}^N} |(D^s u)_\delta(x+th) \cdot h| \cdot |\varphi(x+th) - \varphi(x)| \, dx \\ &= \int_{\mathbb{R}^N} |(D^s u)_\delta(z) \cdot h| \varphi(z) \, dz + |h| \int_{\mathbb{R}^N} |(D^s u)_\delta(z)| \cdot |\varphi(z) - \varphi(z-th)| \, dz.\end{aligned}$$

By Fubini's theorem, we deduce that

$$\begin{aligned}&\int_{\mathbb{R}^N} \int_{\mathbb{R}^N} \left( \int_0^1 \left| (D^s u)_\delta(x+th) \cdot \frac{h}{|h|} \right| dt \right) \rho_\epsilon(h) \, dh \, \varphi(x) \, dx \\ &\leq K_{1,N} \int_{\mathbb{R}^N} |(D^s u)_\delta(z)| \varphi(z) \, dz + \int_0^1 \int_{\mathbb{R}^N} \left( \int_{\mathbb{R}^N} |\varphi(z) - \varphi(z-th)| \rho_\epsilon(h) \, dh \right) |(D^s u)_\delta(z)| \, dz \, dt.\end{aligned}$$

In the second term we can replace the variable $h$ by $-h$ and switch the letter $z$ to $x$. Hence,

$$\begin{aligned}&\int_{\mathbb{R}^N} \left( \int_{\mathbb{R}^N} \frac{|u_\delta(x+h) - u_\delta(x) - (\nabla u)_\delta(x) \cdot h|}{|h|} \rho_\epsilon(h) \, dh \right) \varphi(x) \, dx \\ &\leq K_{1,N} \int_{\mathbb{R}^N} |(D^s u)_\delta(x)| \varphi(x) \, dx \\ &+ \int_0^1 \int_{\mathbb{R}^N} \left( \int_{\mathbb{R}^N} |\varphi(x+th) - \varphi(x)| \rho_\epsilon(h) \, dh \right) |(D^s u)_\delta(x)| \, dx \, dt \\ &+ \|\varphi\|_{L^\infty(\mathbb{R}^N)} \int_0^1 \int_{\mathbb{R}^N} \left( \int_{\mathbb{R}^N} |(\nabla u)_\delta(x+th) - (\nabla u)_\delta(x)| \, dx \right) \rho_\epsilon(h) \, dh \, dt.\end{aligned}$$

Since $|(D^s u)_\delta| \leq |D^s u| * \psi_\delta$, letting $\delta$ tend to zero we deduce (2.1). $\square$

*Proof of Proposition 2.1.* We may assume that $\varphi$ is non-negative. Since $\varphi$ is bounded, there exists $M \geq 0$ such that $\varphi \leq M$ in $\mathbb{R}^N$. We now write

$$\int_{\mathbb{R}^N} R_\epsilon^1 u \, \varphi = M \int_{\mathbb{R}^N} R_\epsilon^1 u - \int_{\mathbb{R}^N} R_\epsilon^1 u \, (M - \varphi).$$

By Lemmas 2.2 and 2.3, we get

$$\liminf_{\epsilon \to 0} \int_{\mathbb{R}^N} R^1_\epsilon u\, \varphi \geq M K_{1,N} |D^s u|(\mathbb{R}^N) - K_{1,N} \int_{\mathbb{R}^N} (M - \varphi)\, \mathrm{d}|D^s u| = K_{1,N} \int_{\mathbb{R}^N} \varphi\, \mathrm{d}|D^s u|.$$

Since the reverse inequality holds for the limsup, the conclusion follows. □

**Corollary 2.4.** *Let $u \in BV(\mathbb{R}^N)$. For every open set $W \subset \mathbb{R}^N$, we have*

$$K_{1,N} |D^s u|(W) \leq \liminf_{\epsilon \to 0} \int_W R^1_\epsilon u \leq \limsup_{\epsilon \to 0} \int_W R^1_\epsilon u \leq K_{1,N} |D^s u|(\overline{W}).$$

*Proof.* We prove the last inequality; the first one can be proved along the same idea using Proposition 2.1 instead of Lemma 2.3. For this purpose, take a sequence of continuous functions $(\varphi_n)_{n \in \mathbb{N}}$ such that $0 \leq \varphi_n \leq 1$ in $\mathbb{R}^N$, $\varphi_n = 1$ on $\overline{W}$ and converging pointwisely to 0 in $\mathbb{R}^N \setminus \overline{W}$. For every $\epsilon > 0$ and $n \in \mathbb{N}$, we have

$$\int_W R^1_\epsilon u \leq \int_{\mathbb{R}^N} R^1_\epsilon u\, \varphi_n.$$

Letting $\epsilon$ tend to zero, by Lemma 2.3 we have that

$$\limsup_{\epsilon \to 0} \int_W R^1_\epsilon u \leq K_{1,N} \int_{\mathbb{R}^N} \varphi_n\, \mathrm{d}|D^s u|.$$

The conclusion follows from the Dominated Convergence Theorem as $n$ tends to infinity. □

*Proof of Theorem 1.4.* The case $\Omega = \mathbb{R}^N$ follows from Proposition 2.1 applied to the test function $\varphi = 1$. We may thus assume that $\Omega$ is an open, bounded and smooth subset of $\mathbb{R}^N$. Given $u \in BV(\Omega)$, we may extend $u$ as a function in $\mathbb{R}^N$, still denoted by $u$, such that $u \in BV(\mathbb{R}^N)$ and $|D^s u|(\partial \Omega) = 0$ [3, Proposition 3.21]. We now consider a function $f_\epsilon : \mathbb{R}^N \times \mathbb{R}^N \to \mathbb{R}$ such that, for every $x \neq y$,

$$f_\epsilon(x, y) = \frac{|u(x) - u(y) - \nabla u(x) \cdot (x - y)|}{|x - y|} \rho_\epsilon(x - y).$$

By an affine change of variables, we have

$$R^1_\epsilon u(x) = \int_{\mathbb{R}^N} f_\epsilon(x, y)\, \mathrm{d}y.$$

Since $|D^s u|(\partial \Omega) = 0$, it thus follows from Corollary 2.4 that

$$\limsup_{\epsilon \to 0} \int_\Omega \int_\Omega f_\epsilon \leq \lim_{\epsilon \to 0} \int_\Omega \int_{\mathbb{R}^N} f_\epsilon = K_{1,N} |D^s u|(\Omega).$$

To obtain the reverse inequality for the liminf, take an open set $V \supseteq \Omega$. We estimate

$$\begin{aligned}
\int_\Omega \int_{\mathbb{R}^N} f_\epsilon &= \int_\Omega \int_\Omega f_\epsilon + \int_\Omega \int_{V \setminus \Omega} f_\epsilon + \int_\Omega \int_{\mathbb{R}^N \setminus V} f_\epsilon \\
&\leq \int_\Omega \int_\Omega f_\epsilon + \int_{V \setminus \overline{\Omega}} \int_{\mathbb{R}^N} f_\epsilon + \int_\Omega \int_{\mathbb{R}^N \setminus V} f_\epsilon.
\end{aligned} \tag{2.2}$$

Taking $r > 0$ such that $d(\mathbb{R}^N \setminus V, \Omega) > r$, we have

$$\int_\Omega \int_{\mathbb{R}^N \setminus V} f_\epsilon \leq \left( \frac{2}{r} \|u\|_{L^1(\mathbb{R}^N)} + \|\nabla u\|_{L^1(\mathbb{R}^N)} \right) \int_{\mathbb{R}^N \setminus B_r(0)} \rho_\epsilon.$$

Thus, as $\epsilon$ tends to zero in estimate (2.2), by Corollary 2.4 we get

$$K_{1,N} |D^s u|(\Omega) \leq \liminf_{\epsilon \to 0} \int_\Omega \int_\Omega f_\epsilon + K_{1,N} |D^s u|(V \setminus \overline{\Omega}).$$

Minimizing the right-hand side with respect to $V$, it follows from the Monotone Set Lemma that
$$K_{1,N}|D^s u|(\Omega) \leq \liminf_{\epsilon \to 0} \int_\Omega \int_\Omega f_\epsilon. \qquad \square$$

## 3. Proofs of Theorem 1.1 and Corollary 1.2

Given $u \in BV(\mathbb{R}^N)$, define the function $S_\epsilon u : \mathbb{R}^N \to \mathbb{R}$ by

(3.1) $$S_\epsilon u(x) = \left| \int_{\mathbb{R}^N} \frac{u(x+h) - u(x)}{|h|} \rho_\epsilon(h) \, dh \right|.$$

Using a regularization argument and the Fundamental Theorem of Calculus (as in the proof of Lemma 2.1 in [23], for example), one can show that

(3.2) $$\|S_\epsilon u\|_{L^1(\mathbb{R}^N)} \leq \int_{\mathbb{R}^N} \int_{\mathbb{R}^N} \frac{|u(x+h) - u(x)|}{|h|} \rho_\epsilon(h) \, dh \, dx \leq |Du|(\mathbb{R}^N),$$

so that $S_\epsilon u$ is well-defined and uniformly bounded in $L^1(\mathbb{R}^N)$ with respect to $\epsilon$. Moreover, as $\rho_\epsilon$ is even, the function $h \mapsto \frac{h}{|h|} \rho_\epsilon(h)$ is odd, and so we additionally have the pointwise estimate

(3.3) $$0 \leq S_\epsilon u \leq R^1_\epsilon u.$$

Let us recall that, for any $u \in BV(\Omega)$, the singular part $D^s u$ of the distributional derivative satisfies the property that, for every measurable subset $E \subset \mathbb{R}^N$ such that $\mathcal{H}^{N-1}(E) = 0$, one has $|D^s u|(E) = 0$. We thus have a Lebesgue decomposition of the measure $D^s u$ of the form
$$D^s u = D^j u + D^c u,$$
where the jump part $D^j u$ is carried by a measurable set which is $\sigma$-finite with respect to the Hausdorff measure $\mathcal{H}^{N-1}$, and $D^c u$ is the Cantor part. We say that $u$ belongs to the class of functions of special bounded variation $SBV(\Omega)$ when $D^c u = 0$.

This jump–Cantor terminology comes from the fact that, by the Federer-Vol'pert decomposition of $D^j u$ [3, Theorem 3.78], there exists an $(N-1)$-dimensional rectifiable set $J_u$ such that there exists a triple $(u^+, u^-, \nu_u) \in \mathbb{R} \times \mathbb{R} \times \mathbb{S}^{N-1}$ of Borel functions defined on $J_u$ which satisfies

(3.4) $$D^j u = (u^+ - u^-) \otimes \nu_u \, \mathcal{H}^{N-1} \llcorner J_u.$$

In dimension one, $D^j u$ is a countable combination of Dirac masses, each one indicating a jump discontinuity of the function. The distributional derivative of the Cantor function contains only the Cantor part, and is supported on the standard Cantor set.

In what follows, we need some additional properties on the jump set $J_u$. Since $J_u$ is rectifiable, it is contained in a countable union of graphs of Lipschitz functions, and for $\mathcal{H}^{N-1}$-almost every $x_0 \in J_u$, the following properties hold:

(i) [21, Theorem 10.2] there exists an approximate tangent plane $T_{x_0}(J_u)$ at $x_0$ and
$$\lim_{r \to 0} \frac{\mathcal{H}^{N-1}(J_u \cap B_r(x_0))}{\omega_{N-1} r^{N-1}} = 1,$$
and [3, Theorem 3.78] $\nu_u(x_0)$ is orthogonal to $T_{x_0}(J_u)$;

(ii) [21, Corollary 6.5 and Proposition 10.5] up to a rotation, there exists a Lipschitz function $\gamma : \mathbb{R}^{N-1} \to \mathbb{R}$ such that $x_0$ is contained in the graph $\operatorname{gr} \gamma$, the approximate tangent planes $T_{x_0}(\operatorname{gr} \gamma)$ and $T_{x_0}(J_u)$ exist and coincide,
$$\lim_{r \to 0} \frac{\mathcal{H}^{N-1}(\operatorname{gr} \gamma \cap B_r(x_0))}{\omega_{N-1} r^{N-1}} = 1,$$

88 AUGUSTO C. PONCE AND DANIEL SPECTOR

and
$$\lim_{r \to 0} \frac{\mathcal{H}^{N-1}((\operatorname{gr} \gamma \triangle J_u) \cap B_r(x_0))}{r^{N-1}} = 0,$$
where $A \triangle B = (A \setminus B) \cup (B \setminus A)$ denotes the symmetric difference between the sets $A$ and $B$;

(iii) [21, Theorem 5.16] $x_0$ is a Lebesgue point of $u^+$ and $u^-$ with respect to the measure $\mathcal{H}^{N-1} \lfloor_{J_u}$, and the precise representatives are $u^+(x_0)$ and $u^-(x_0)$, respectively.

Given $\kappa > 0$, denote by $T_\kappa : \mathbb{R} \to \mathbb{R}$ the truncation function at levels $\pm \kappa$:
$$T_\kappa(s) = \begin{cases} \kappa & \text{if } s > \kappa, \\ s & \text{if } -\kappa \leq s \leq \kappa, \\ -\kappa & \text{if } s < -\kappa. \end{cases}$$

For every $u \in BV(\mathbb{R}^N)$, by Lipschitz continuity of $T_\kappa$ one deduces using a smooth approximation of $u$ that $T_\kappa(u) \in BV(\mathbb{R}^N)$ and $|D(T_\kappa(u))|(\mathbb{R}^N) \leq |Du|(\mathbb{R}^N)$. It follows from the Federer-Vol'pert chain rule formula for $BV$ functions [2, Corollary 3.1; 3, Theorem 3.96] that $J_{T_\kappa(u)} \subset J_u$ and

(3.5) $$D^j(T_\kappa(u)) = (T_\kappa(u^+) - T_\kappa(u^-)) \otimes \nu_u \, \mathcal{H}^{N-1} \lfloor_{J_u}.$$

From the Dominated Convergence Theorem, we then have that

(3.6) $$\lim_{k \to \infty} |D^j(T_\kappa(u)) - D^j u|(\mathbb{R}^N) = 0.$$

We now have the ingredients to establish

**Lemma 3.1.** *Let $u \in BV(\mathbb{R}^N)$. For every non-negative continuous function with compact support $\varphi : \mathbb{R}^N \to \mathbb{R}$, we have*
$$K_{1,N} \int_{\mathbb{R}^N} \varphi \, \mathrm{d}|D^j u| \leq \liminf_{\epsilon \to 0} \int_{\mathbb{R}^N} S_\epsilon u \, \varphi.$$

*Proof.* Since $u \in BV(\mathbb{R}^N)$, the inequality (3.2) implies that the family $(S_\epsilon u)_{\epsilon > 0}$ is bounded in $(C_0(\mathbb{R}^N))'$, where $C_0(\mathbb{R}^N)$ denotes the Banach space of continuous real functions in $\mathbb{R}^N$ that converge uniformly to zero at infinity. Given $\varphi$, take a sequence $(\epsilon_n)_{n \in \mathbb{N}}$ converging to zero such that
$$\lim_{n \to \infty} \int_{\mathbb{R}^N} S_{\epsilon_n} u \, \varphi = \liminf_{\epsilon \to 0} \int_{\mathbb{R}^N} S_\epsilon u \, \varphi$$
and $S_{\epsilon_n} u \overset{*}{\rightharpoonup} \lambda$ for some finite Borel measure $\lambda$ in $\mathbb{R}^N$. The existence of $\lambda$ follows from the Riesz Representation Theorem [21, Theorem 4.7].

We prove that $\lambda \geq K_{1,N} |D^j u|$. Given $x_0 \in J_u$ that satisfies Assertions (i)–(iii) and is contained in the graph $\operatorname{gr} \gamma$ of a Lipschitz function as above, we partition the space $\mathbb{R}^N$ in two parts:
$$\mathbb{R}_+^N = \{(x,t) \in \mathbb{R}^{N-1} \times \mathbb{R} : t > \gamma(x)\} \quad \text{and} \quad \mathbb{R}_-^N = \{(x,t) \in \mathbb{R}^{N-1} \times \mathbb{R} : t < \gamma(x)\}.$$

Let $g : \mathbb{R}^N \to \mathbb{R}$ be the function defined by
$$g(x) = \begin{cases} u^+(x_0) & \text{if } x \in \mathbb{R}_+^N, \\ u^-(x_0) & \text{if } x \in \mathbb{R}_-^N. \end{cases}$$

We have that $g \in BV_{\mathrm{loc}}(\mathbb{R}^N)$ and
$$Dg = (u^+(x_0) - u^-(x_0)) \otimes \nu_\gamma \, \mathcal{H}^{N-1} \lfloor_{\operatorname{gr} \gamma},$$
where $\nu_\gamma$ is the unit normal vector with respect to the graph of $\gamma$. Even though $g$ does not belong to $BV(\mathbb{R}^N)$, nor to $L^1(\mathbb{R}^N)$, we may define $S_\epsilon g$ and $S_\epsilon(u - g)$ as in the formula (3.1) and they both belong to $L^1_{\text{loc}}(\mathbb{R}^N)$. To see this, let $B_R(0)$ be the ball centered at zero with radius $R > 0$. Then one has the estimate

$$\int_{B_R(0)} |S_\epsilon g(x)| \, \mathrm{d}x \leq \int_{B_R(0)} \int_{B_1(0)} \frac{|g(x+h) - g(x)|}{|h|} \rho_\epsilon(h) \, \mathrm{d}h \, \mathrm{d}x$$
$$+ \int_{B_R(0)} \int_{\mathbb{R}^N \setminus B_1(0)} 2\|g\|_{L^\infty(\mathbb{R}^N)} \rho_\epsilon(h) \, \mathrm{d}h \, \mathrm{d}x$$
$$\leq |Dg|(B_{R+1}(0)) + 2\omega_N R^N \|g\|_{L^\infty(\mathbb{R}^N)}.$$

Thus, as $S_\epsilon u, S_\epsilon g \in L^1_{\text{loc}}(\mathbb{R}^N)$ so also does $S_\epsilon(u - g)$.

By the triangle inequality, we have

(3.7) $$S_\epsilon u \geq S_\epsilon g - S_\epsilon(u - g).$$

Applying the pointwise estimate (3.3) and Lemma 2.3 to $u - g$, we have that

(3.8) $$\limsup_{\epsilon \to 0} \int_{\mathbb{R}^N} S_\epsilon(u-g)\varphi \leq \int_{\mathbb{R}^N} \varphi \, \mathrm{d}|D^s(u-g)|.$$

Moreover, since the image of $g$ has only two elements, by a straightforward computation in the spirit of (1.5) we have the identity

$$\int_{\mathbb{R}^N} S_\epsilon g \, \varphi = \int_{\mathbb{R}^N} \int_{\mathbb{R}^N} \frac{|g(x) - g(y)|}{|x - y|} \rho_\epsilon(x - y) \varphi(x) \, \mathrm{d}y \, \mathrm{d}x.$$

Hence, by the local version of Dávila's result [14, Lemma 2] we get

(3.9) $$\lim_{\epsilon \to 0} \int_{\mathbb{R}^N} S_\epsilon g \, \varphi = K_{1,N} \int_{\mathbb{R}^N} \varphi \, \mathrm{d}|Dg|.$$

Recalling that $S_{\epsilon_n} u \overset{*}{\rightharpoonup} \lambda$ and combining equations (3.7)–(3.9), we deduce that

$$\int_{\mathbb{R}^N} \varphi \, \mathrm{d}\lambda \geq K_{1,N} \int_{\mathbb{R}^N} \varphi \, \mathrm{d}|Dg| - K_{1,N} \int_{\mathbb{R}^N} \varphi \, \mathrm{d}|D^s(u-g)|.$$

Hence, on every compact subset of $\mathbb{R}^N$, we have

$$\lambda \geq K_{1,N}|Dg| - K_{1,N}|D^s(u-g)|,$$

and then, by inner regularity of finite Borel measures, this inequality holds on every Borel subset of $\mathbb{R}^N$.

Denote by $\lambda^j$ the measure $\lambda \llcorner J_u$. Then restricting the previous inequality to $J_u$ and using the Federer-Vol'pert decomposition formula (3.4), we have the following relation between measures:

$$\frac{1}{K_{1,N}} \lambda^j \geq |D^j g| - |D^j(u-g)|$$
$$= |u^+(x_0) - u^-(x_0)| \mathcal{H}^{N-1} \llcorner \text{gr } \gamma$$
$$- |(u^+ - u^+(x_0)) - (u^- - u^-(x_0))| \mathcal{H}^{N-1} \llcorner \text{gr } \gamma \cap J_u$$
$$- |D^j u| \llcorner_{J_u \setminus \text{gr } \gamma} - |u^+(x_0) - u^-(x_0)| \mathcal{H}^{N-1} \llcorner \text{gr } \gamma \setminus J_u.$$

For any given $\kappa > 0$, by the triangle inequality and the structure of the measure $D^j T_\kappa(u)$ given by the chain rule (3.5), we have

$$|D^j u| \llcorner_{J_u \setminus \text{gr } \gamma} \leq |D^j(T_\kappa(u))| \llcorner_{J_u \setminus \text{gr } \gamma} + |D^j(T_\kappa(u)) - D^j u| \llcorner_{J_u \setminus \text{gr } \gamma}$$
$$\leq 2\kappa \mathcal{H}^{N-1} \llcorner_{J_u \setminus \text{gr } \gamma} + |D^j(T_\kappa(u)) - D^j u| \llcorner_{J_u}.$$

The last measure in the right-hand side does not depend on $x_0$. Denoting

$$\mu_\kappa = \frac{1}{K_{1,N}} \lambda^j + |D^j(T_\kappa(u)) - D^j u| \llcorner_{J_u},$$

we thus have that

$$\mu_\kappa \geq |u^+(x_0) - u^-(x_0)|\mathcal{H}^{N-1}\lfloor_{\operatorname{gr}\gamma}$$
$$- |(u^+ - u^+(x_0)) - (u^- - u^-(x_0))|\mathcal{H}^{N-1}\lfloor_{\operatorname{gr}\gamma \cap J_u}$$
$$- 2\kappa\mathcal{H}^{N-1}\lfloor_{J_u\setminus\operatorname{gr}\gamma} - |u^+(x_0) - u^-(x_0)|\mathcal{H}^{N-1}\lfloor_{\operatorname{gr}\gamma\setminus J_u}.$$

By Properties $(i)$–$(iii)$ satisfied by the point $x_0$, we deduce that

$$\liminf_{r\to 0} \frac{\mu_\kappa(B_r(x_0))}{\omega_{N-1}r^{N-1}} \geq |u^+(x_0) - u^-(x_0)|.$$

Since this relation holds for $\mathcal{H}^{N-1}$-almost every point $x_0 \in J_u$, it follows from the Besicovitch differentiation theorem [21, Theorem 6.4] that

$$\frac{1}{K_{1,N}}\lambda^j + |D^j(T_\kappa(u)) - D^j u|\lfloor_{J_u} = \mu_\kappa \geq |u^+ - u^-|\mathcal{H}^{N-1}\lfloor_{J_u}.$$

As $\kappa$ tends to infinity, we deduce from the strong convergence (3.6) that

$$\frac{1}{K_{1,N}}\lambda^j \geq |u^+ - u^-|\mathcal{H}^{N-1}\lfloor_{J_u} = |D^j u|.$$

Therefore,

$$\liminf_{\epsilon\to 0}\int_{\mathbb{R}^N} S_\epsilon u\, \varphi = \int_{\mathbb{R}^N}\varphi\,\mathrm{d}\lambda \geq \int_{\mathbb{R}^N}\varphi\,\mathrm{d}\lambda^j \geq K_{1,N}\int_{\mathbb{R}^N}\varphi\,\mathrm{d}|D^j u|. \qquad \square$$

*Proof of Theorem 1.1.* Let $u \in SBV(\mathbb{R}^N)$. We first prove that, for every bounded continuous function $\varphi : \mathbb{R}^N \to \mathbb{R}$ (not necessarily with compact support), we have

$$(3.10) \qquad K_{1,N}\int_{\mathbb{R}^N}\varphi\,\mathrm{d}|D^j u| \leq \liminf_{\epsilon\to 0}\int_{\mathbb{R}^N} S_\epsilon u\,\varphi.$$

Taking this estimate from granted, we may combine it with the pointwise inequality (3.3) and Lemma 2.3 to get

$$K_{1,N}\int_{\mathbb{R}^N}\varphi\,\mathrm{d}|D^j u| \leq \liminf_{\epsilon\to 0}\int_{\mathbb{R}^N} S_\epsilon u\,\varphi \leq \limsup_{\epsilon\to 0}\int_{\mathbb{R}^N} S_\epsilon u\,\varphi \leq K_{1,N}\int_{\mathbb{R}^N}\varphi\,\mathrm{d}|D^s u|.$$

Since $D^s u = D^j u$, we deduce that

$$(3.11) \qquad \lim_{\epsilon\to 0}\int_{\mathbb{R}^N} S_\epsilon u\,\varphi = K_{1,N}\int_{\mathbb{R}^N}\varphi\,\mathrm{d}|D^s u|.$$

Choosing $\varphi = 1$, we have the conclusion when $\Omega = \mathbb{R}^N$. In the case when $\Omega$ is an open, bounded, smooth subset of $\mathbb{R}^N$, we can proceed by extension to $\mathbb{R}^N$ along the lines of the proof of Theorem 1.4.

To prove (3.10), it suffices to consider the case of a non-negative function $\varphi$. Given a sequence of continuous functions with compact support $(\psi_n)_{n\in\mathbb{N}}$ converging pointwise to 1 and such that $0 \leq \psi_n \leq 1$ in $\mathbb{R}^N$, we have

$$\int_{\mathbb{R}^N} S_\epsilon u\,\varphi\psi_n \leq \int_{\mathbb{R}^N} S_\epsilon u\,\varphi.$$

Letting $\epsilon$ tend to zero, we deduce from Lemma 3.1 that

$$K_{1,N}\int_{\mathbb{R}^N}\varphi\psi_n\,\mathrm{d}|D^j u| \leq \liminf_{\epsilon\to 0}\int_{\mathbb{R}^N} S_\epsilon u\,\varphi.$$

The conclusion follows from the Dominated Convergence Theorem as $n$ tends to infinity. $\square$

In the previous proof, we have established the strict convergence of $(S_\epsilon u)_{\epsilon>0}$. This conclusion still holds under a weaker assumption on the family $(\rho_\epsilon)_{\epsilon>0}$:

**Proposition 3.2.** *Let $(\rho_\epsilon)_{\epsilon>0} \subset L^1_{\mathrm{loc}}(\mathbb{R}^N)$ be a family of non-negative radial functions in $\mathbb{R}^N$ that satisfies* (1.7) *and*

$$(3.12) \qquad \lim_{\epsilon \to 0} \int_{|h|>\delta} \frac{\rho_\epsilon(h)}{|h|}\,\mathrm{d}h = 0,$$

*for every $\delta > 0$. Then, for every $u \in SBV(\mathbb{R}^N)$ and every bounded continuous function $\varphi: \mathbb{R}^N \to \mathbb{R}$, we have*

$$\lim_{\epsilon \to 0} \int_{\mathbb{R}^N} S_\epsilon u\,\varphi = K_{1,N} \int_{\mathbb{R}^N} \varphi\,\mathrm{d}|D^s u|.$$

*Proof.* The function $S_\epsilon u$ is still well-defined in this case since $u \in L^1(\mathbb{R}^N)$. Given $R > 0$, let

$$S_{\epsilon,R} u(x) = \left| \int_{\mathbb{R}^N} \frac{u(x+h) - u(x)}{|h|} \rho_\epsilon(h) \chi_{B_R(0)}(h)\,\mathrm{d}h \right|.$$

By the triangle inequality, for every $x \in \mathbb{R}^N$ we have

$$|S_\epsilon u(x) - S_{\epsilon,R} u(x)| \leq \int_{\mathbb{R}^N \setminus B_R(0)} \frac{|u(x+h) - u(x)|}{|h|} \rho_\epsilon(h)\,\mathrm{d}h.$$

By Fubini's theorem, we thus get

$$\|S_\epsilon u - S_{\epsilon,R} u\|_{L^1(\mathbb{R}^N)} \leq 2\|u\|_{L^1(\mathbb{R}^N)} \int_{\mathbb{R}^N \setminus B_R(0)} \frac{\rho_\epsilon(h)}{|h|}\,\mathrm{d}h.$$

In particular, for every bounded continuous function $\varphi: \mathbb{R}^N \to \mathbb{R}$, we have

$$\left| \int_{\mathbb{R}^N} S_\epsilon u\,\varphi - \int_{\mathbb{R}^N} S_{\epsilon,R} u\,\varphi \right| \leq 2\|\varphi\|_{L^\infty(\mathbb{R}^N)} \|u\|_{L^1(\mathbb{R}^N)} \int_{\mathbb{R}^N \setminus B_R(0)} \frac{\rho_\epsilon(h)}{|h|}\,\mathrm{d}h.$$

Letting $\epsilon$ tend to zero, the conclusion follows from formula (3.11) applied to the family $(\rho_\epsilon \chi_{B_R(0)})_{\epsilon > 0}$. $\square$

*Proof of Corollary 1.2.* It suffices to apply Proposition 3.2 to

$$\rho_\epsilon(h) = \frac{\epsilon}{\mathcal{H}^{N-1}(\mathbb{S}^{N-1})} \frac{1}{|h|^{N-\epsilon}},$$

and then take $\epsilon = 1 - \alpha$. Since $c_1 C_{1,N} = 2/\pi$, we have the conclusion. $\square$

## 4. A COMPUTATION INVOLVING THE CANTOR FUNCTION

The goal of this section is to motivate why one should expect a non-trivial contribution from the Cantor part of the derivative of a $BV$ function for the functional in Theorem 1.1. We focus our attention on the standard Cantor function: this is the unique fixed point of the contraction map $T: X \to X$ defined by

$$T(v)(x) = \begin{cases} v(3x)/2 & \text{if } 0 \leq x \leq 1/3, \\ 1/2 & \text{if } 1/3 < x < 2/3, \\ (v(3x-2)+1)/2 & \text{if } 2/3 \leq x \leq 1, \end{cases}$$

where $X$ is the complete metric space formed by all continuous functions $v: [0,1] \to \mathbb{R}$ such that $v(0) = 0$ and $v(1) = 1$, equipped with the uniform distance. Since the Cantor function is non-decreasing, it belongs to $BV(0,1)$.

**Proposition 4.1.** *If $u: [0,1] \to \mathbb{R}$ is the standard Cantor function, then, for every $0 < \epsilon < 1/6$ and every non-negative even function $\rho_\epsilon \in L^1(\mathbb{R})$ with $\mathrm{supp}\,\rho_\epsilon \subset [-\epsilon, \epsilon]$, we have*

$$\int_0^1 \left| \int_0^1 \frac{u(x) - u(y)}{|x-y|} \rho_\epsilon(x-y)\,\mathrm{d}y \right|\mathrm{d}x \geq \frac{1}{36} \int_{3\epsilon/4}^\epsilon \rho_\epsilon.$$

In the case where $\rho_\epsilon = \chi_{(-\epsilon,\epsilon)}/2\epsilon$, the right-hand side is constant and positive, and we deduce that

$$\liminf_{\epsilon \to 0} \int_0^1 \left| \int_0^1 \frac{u(x) - u(y)}{|x-y|} \rho_\epsilon(x-y) \, \mathrm{d}x \right| \mathrm{d}y \geq \frac{1}{36 \cdot 8}.$$

We begin with the following identity based on Fubini's theorem:

**Lemma 4.2.** *Let $u \in BV(0,1) \cap C^0[0,1]$ be a non-decreasing function, and let $\rho \in L^1(\mathbb{R})$ be a non-negative even function with compact support. For every measurable subset $A \subset (0,1)$ such that $A - \operatorname{supp} \rho \subset (0,1)$, we have*

$$\int_A \left( \int_0^1 \frac{|u(x)-u(y)|}{|x-y|} \rho(x-y) \, \mathrm{d}y \right) \mathrm{d}x = \int_0^1 \left( \int_{\mathbb{R}_+} \frac{\left|(z-h, z+h) \cap A\right|}{2h} 2\rho(h) \, \mathrm{d}h \right) \mathrm{d}|Du|(z).$$

*Proof.* We first extend $u$ as a continuous function in $\mathbb{R}$, which we still denote by $u$, by taking $u(x) = u(1)$ for $x > 1$ and $u(x) = u(0)$ for $x < 0$. In particular, the measure $Du$ is supported on $[0,1]$. Since $A - \operatorname{supp} \rho \subset (0,1)$, we have

$$\int_A \left( \int_0^1 \frac{|u(x)-u(y)|}{|x-y|} \rho(x-y) \, \mathrm{d}y \right) \mathrm{d}x = \int_A \left( \int_{\mathbb{R}} \frac{|u(x)-u(y)|}{|x-y|} \rho(x-y) \, \mathrm{d}y \right) \mathrm{d}x.$$

Denote this common quantity by $I(u)$. Making the change of variable $y = x+h$ with respect to $y$, we get

$$I(u) = \int_A \left( \int_{\mathbb{R}} \frac{|u(x+h)-u(x)|}{|h|} \rho(h) \, \mathrm{d}h \right) \mathrm{d}x.$$

Since $u$ is continuous and non-decreasing, the measure $Du$ does not contain atoms and is non-negative. For every $h > 0$, we then have

$$|u(x+h) - u(x)| = u(x+h) - u(x) = \int_x^{x+h} \mathrm{d}Du = |Du|(x, x+h).$$

A similar computation holds for $h < 0$. Since $\rho$ is even, we get

$$I(u) = \int_A \left( \int_{\mathbb{R}_+} \frac{|Du|(x-h, x+h)}{h} \rho(h) \, \mathrm{d}h \right) \mathrm{d}x = \int_{\mathbb{R}_+} \left( \int_A \frac{|Du|(x-h, x+h)}{h} \, \mathrm{d}x \right) \rho(h) \, \mathrm{d}h.$$

We apply again Fubini's theorem to reach the conclusion. Indeed, for every $h > 0$ we have

$$\int_A |Du|(x-h, x+h) \, \mathrm{d}x = \int_A \left( \int_{x-h < z < x+h} \mathrm{d}|Du|(z) \right) \mathrm{d}x$$

$$= \int_{A+(-h,h)} \left( \int_{x \in A,\, z-h < x < z+h} \mathrm{d}x \right) \mathrm{d}|Du|(z)$$

$$= \int_{A+(-h,h)} \left|(z-h, z+h) \cap A\right| \mathrm{d}|Du|(z).$$

Note that for $z \notin A + (-h, h)$, the set $(z-h, z+h) \cap A$ is empty. We thus get

$$\int_A |Du|(x-h, x+h) \, \mathrm{d}x = \int_{\mathbb{R}} \left|(z-h, z+h) \cap A\right| \mathrm{d}|Du|(z) = \int_0^1 \left|(z-h, z+h) \cap A\right| \mathrm{d}|Du|(z).$$

The proof of the lemma is complete. $\square$

*Proof of Proposition 4.1.* Given $0 < \epsilon < 1/6$, let $n \in \mathbb{N}_*$ be the greatest integer such that $2\epsilon \leq 1/3^n$. Taking $A_n$ to be the union of the constant sections of the Cantor function of length $1/3^n$, we have that

$$\int_0^1 \left| \int_0^1 \frac{u(x) - u(y)}{|x-y|} \rho_\epsilon(x-y) \, dy \right| dx \geq \int_{A_n} \left| \int_0^1 \frac{u(x) - u(y)}{|x-y|} \rho_\epsilon(x-y) \, dy \right| dx$$
$$= \int_{A_n} \left( \int_0^1 \frac{|u(x) - u(y)|}{|x-y|} \rho_\epsilon(x-y) \, dy \right) dx.$$

Indeed, the properties $\operatorname{supp} \rho_\epsilon \subset [-\epsilon, \epsilon]$ and $2\epsilon \leq 1/3^n$ ensure that, for each $x \in A_n$, the function

$$y \in (0,1) \longmapsto \frac{u(x) - u(y)}{|x-y|} \rho_\epsilon(x-y)$$

has constant sign. It thus follows from the previous lemma that
(4.1)
$$\int_0^1 \left| \int_0^1 \frac{u(x) - u(y)}{|x-y|} \rho_\epsilon(x-y) \, dy \right| dx \geq \int_0^1 \left( \int_{\mathbb{R}_+} \frac{|(z-h, z+h) \cap A_n|}{h} \rho_\epsilon(h) \, dh \right) d|Du|(z).$$

In the construction of the Cantor set, let $\mathcal{I}_n$ be the collection of closed intervals of length $1/3^n$ that contain the Cantor set. For example, $\mathcal{I}_1 = \{[0, 1/3], [2/3, 1]\}$. Let $\mathcal{J}_{n+2}$ be the subcollection of intervals in $\mathcal{I}_{n+2}$ that intersect the closure of an interval of length $1/3^n$ that is *removed* during the construction of the Cantor set. For example, $\mathcal{J}_3 = \{[1/3 - 1/27, 1/3], [2/3, 2/3 + 1/27]\}$, and the number of elements of $\mathcal{J}_{n+2}$ is $2^n$.

Given $J \in \mathcal{J}_{n+2}$ and $z \in J$, we have that the set $(z-h, z+h) \cap A_n$ is non-empty for every $h > 1/3^{n+2}$. In particular, for every $1/3^{n+2} \leq h \leq 1/3^n$ we have

$$\frac{|(z-h, z+h) \cap A_n|}{h} \geq \frac{h - 1/3^{n+2}}{h} = 1 - \frac{1/3^{n+2}}{h}.$$

Since $1/3^{n+1} < 2\epsilon \leq 1/3^n$, for every $3\epsilon/4 \leq h \leq \epsilon$ we deduce that

$$\frac{|(z-h, z+h) \cap A_n|}{h} \geq 1 - \frac{2\epsilon/3}{3\epsilon/4} = \frac{1}{9}.$$

Hence, for every $z \in J$ we have

$$\int_{\mathbb{R}_+} \frac{|(z-h, z+h) \cap A_n|}{h} \rho_\epsilon(h) \, dh \geq \int_{3\epsilon/4}^\epsilon \frac{|(z-h, z+h) \cap A_n|}{h} \rho_\epsilon(h) \, dh \geq \frac{1}{9} \int_{3\epsilon/4}^\epsilon \rho_\epsilon.$$

Writing $J = [a_i, b_i] \in \mathcal{J}_{n+2}$, we also have

$$|Du|(J) = u(b_i) - u(a_i) = \frac{1}{2^{n+2}}.$$

We then deduce that

$$\int_J \left( \int_{\mathbb{R}_+} \frac{|(z-h, z+h) \cap A_n|}{h} \rho_\epsilon(h) \, dh \right) d|Du|(z) \geq \frac{1}{2^{n+2}} \cdot \frac{1}{9} \int_{3\epsilon/4}^\epsilon \rho_\epsilon.$$

Since the family $\mathcal{J}_{n+2}$ has $2^n$ elements, using estimate (4.1) we get

$$\int_0^1 \left| \int_0^1 \frac{u(x) - u(y)}{|x-y|} \rho_\epsilon(x-y) \, dy \right| dx$$
$$\geq \sum_{J \in \mathcal{J}_{n+2}} \int_J \left( \int_{\mathbb{R}_+} \frac{|(z-h, z+h) \cap A_n|}{h} \rho_\epsilon(h) \, dh \right) d|Du|(z) \geq \frac{1}{36} \int_{3\epsilon/4}^\epsilon \rho_\epsilon.$$

This gives the conclusion. □

## 5. Connection with a BMO-type norm

Inspired from [8], in the paper [1, Section 4.3] the authors associate to $u \in L^1_{\text{loc}}(\mathbb{R}^N)$ a semi-norm

$$\kappa_\epsilon(u) = \epsilon^{N-1} \sup_{\mathcal{Q}_\epsilon} \sum_{Q' \in \mathcal{Q}_\epsilon} \fint_{Q'} \left| u(x) - \fint_{Q'} u(y) \, \mathrm{d}y \right| \mathrm{d}x,$$

where the supremum is taken over all families $\mathcal{Q}_\epsilon$ of disjoint cubes with side length $\epsilon$, with arbitrary orientation and cardinality. In particular, for $u = \chi_A$, where $A \subset \mathbb{R}^N$ is a subset of finite perimeter, they show that one has the convergence [1, Eq. (4.4)]

$$\lim_{\epsilon \to 0} \kappa_\epsilon(\chi_A) = \frac{1}{2} |D\chi_A|(\mathbb{R}^N).$$

This result has been extended by Fusco, Moscariello and Sbordone [18, Theorem 3.3] to $u \in SBV(\mathbb{R}^N)$ such that either $\nabla u \equiv 0$ or $|\overline{J}_u| = 0$, where they prove that for such $u$ one has

$$\lim_{\epsilon \to 0} \kappa_\epsilon(u) = \frac{1}{4} \int_{\mathbb{R}^N} |\nabla u| + \frac{1}{2} |D^s u|(\mathbb{R}^N).$$

It is not difficult to see that any limit of the family $(\kappa_\epsilon(u))_{\epsilon > 0}$ is comparable to the total variation of $Du$:

**Proposition 5.1.** *For every $u \in BV(\mathbb{R}^N)$, we have*

$$\frac{1}{4}|Du|(\mathbb{R}^N) \leq \liminf_{\epsilon \to 0} \kappa_\epsilon(u) \leq \limsup_{\epsilon \to 0} \kappa_\epsilon(u) \leq \frac{1}{2}|Du|(\mathbb{R}^N).$$

*Proof.* For the upper bound one can apply the Poincaré inequality for functions of bounded variation to obtain

$$\sum_{Q' \in \mathcal{Q}_\epsilon} \epsilon^{N-1} \fint_{Q'} \left| u(x) - \fint_{Q'} u(y) \, \mathrm{d}y \right| \mathrm{d}x \leq \sum_{Q' \in \mathcal{Q}_\epsilon} \frac{1}{2} |Du|(Q') \leq \frac{1}{2}|Du|(\mathbb{R}^N).$$

For the lower bound, mollification lowers the energy following an observation by E. Stein, as for the functionals introduced by Bourgain, Brezis and Mironescu [7]; cf. [10]. Indeed, letting $u_\delta = u * \psi_\delta$, by Fubini's theorem one has

$$\fint_{Q'} \left| u_\delta(x) - \fint_{Q'} u_\delta(y) \, \mathrm{d}y \right| \mathrm{d}x \leq \int_{\mathbb{R}^N} \psi_\delta(z) \left( \fint_{Q'-z} \left| u(x) - \fint_{Q'-z} u(y) \, \mathrm{d}y \right| \mathrm{d}x \right) \mathrm{d}z,$$

and this implies that

$$\sum_{Q' \in \mathcal{Q}_\epsilon} \epsilon^{N-1} \fint_{Q'} \left| u_\delta(x) - \fint_{Q'} u_\delta(y) \, \mathrm{d}y \right| \mathrm{d}x \leq \kappa_\epsilon(u) \int_{\mathbb{R}^N} \psi_\delta(z) \, \mathrm{d}z = \kappa_\epsilon(u).$$

Hence,

$$\kappa_\epsilon(u_\delta) \leq \kappa_\epsilon(u),$$

so that first applying the result of Fusco, Moscariello and Sbordone [18, Theorem 3.3] to the $W^{1,1}(\mathbb{R}^N)$ (smooth) function $u_\delta$ and then using the lower semicontinuity of the total variation with respect to the strict convergence of the family $(u_\delta)_{\delta > 0}$ in $BV(\mathbb{R}^N)$, one deduces the lower bound

$$\frac{1}{4}|Du|(\mathbb{R}^N) \leq \liminf_{\epsilon \to 0} \kappa_\epsilon(u). \qquad \square$$

We can also show that, up to an affine change in $\epsilon$, the energies $\kappa_\epsilon$ are comparable at every scale to a special case of those originally introduced by Bourgain, Brezis and Mironescu [7]. To this end, we first introduce an equivalent energy to $\kappa_\epsilon$, the functional

$$\kappa'_\epsilon(u) = \epsilon^{N-1} \sup_{\mathcal{Q}_\epsilon} \sum_{Q' \in \mathcal{Q}_\epsilon} \fint_{Q'} \fint_{Q'} |u(x) - u(y)|\, dy\, dx.$$

As is the case with semi-norms of functions with bounded mean oscillation (BMO) [13], one has

$$\kappa_\epsilon(u) \leq \kappa'_\epsilon(u) \leq 2\kappa_\epsilon(u),$$

and so it suffices to establish the equivalence of $\kappa'_\epsilon$. We now make the choice of the family of mollifiers $(\rho_\epsilon)_{\epsilon>0}$,

(5.1) $$\rho_\epsilon(h) = \frac{C_N}{\epsilon^{N+1}} |h| \chi_{B_\epsilon(0)}(h),$$

which gives
(5.2)
$$\int_{\mathbb{R}^N} \int_{\mathbb{R}^N} \frac{|u(x) - u(y)|}{|x-y|} \rho_\epsilon(x-y)\, dy\, dx = \frac{C_N \omega_N}{\epsilon} \int_{\mathbb{R}^N} \left( \fint_{B_\epsilon(x)} |u(x) - u(y)|\, dy \right) dx.$$

Then, the equivalence we assert is

**Proposition 5.2.** *Let $u \in L^1_{\mathrm{loc}}(\mathbb{R}^N)$. For every $\epsilon > 0$, we have*

$$C' \kappa'_{\frac{\epsilon}{\sqrt{N}}}(u) \leq \frac{1}{\epsilon} \int_{\mathbb{R}^N} \left( \fint_{B_\epsilon(x)} |u(x) - u(y)|\, dy \right) dx \leq C'' \kappa'_{3\epsilon}(u),$$

*for some constants $C', C'' > 0$ depending on $N$.*

*Proof.* Given a cube $Q'$ with side length $\frac{\epsilon}{\sqrt{N}}$, for every $x \in Q'$ we have that $Q' \subset B_\epsilon(x)$. Hence,

$$\int_{Q'} \left( \int_{Q'} |u(x) - u(y)|\, dy \right) dx \leq \int_{Q'} \left( \int_{B_\epsilon(x)} |u(x) - u(y)|\, dy \right) dx.$$

It thus follows that

$$\left( \frac{\epsilon}{\sqrt{N}} \right)^{N-1} \sum_{Q' \in \mathcal{Q}_{\frac{\epsilon}{\sqrt{N}}}} \fint_{Q'} \fint_{Q'} |u(x) - u(y)|\, dy\, dx \leq \frac{C_1}{\epsilon} \sum_{Q' \in \mathcal{Q}_{\frac{\epsilon}{\sqrt{N}}}} \int_{Q'} \left( \fint_{B_\epsilon(x)} |u(x) - u(y)|\, dy \right) dx$$

$$\leq \frac{C_1}{\epsilon} \int_{\mathbb{R}^N} \left( \fint_{B_\epsilon(x)} |u(x) - u(y)|\, dy \right) dx.$$

Taking the supremum with respect to all families of cubes $\mathcal{Q}_{\frac{\epsilon}{\sqrt{N}}}$, we deduce the first inequality.

To prove the second inequality, write $\mathbb{R}^N$ as a union of disjoint cubes $\mathcal{Q}_{3\epsilon}$ with side length $3\epsilon$. Divide each cube $Q' \in \mathcal{Q}_{3\epsilon}$ in $3^N$ cubes with disjoint interiors and side length $\epsilon$. Denote by $Q''$ the inner cube; this is the only cube that does not intersect $\partial Q'$. For every $x \in Q''$, we have $B_\epsilon(x) \subset Q'$. Hence,

$$\int_{Q''} \left( \int_{B_\epsilon(x)} |u(x) - u(y)|\, dy \right) dx \leq \int_{Q'} \left( \int_{Q'} |u(x) - u(y)|\, dy \right) dx,$$

which implies that

$$\frac{1}{\epsilon} \sum_{Q' \in \mathcal{Q}_{3\epsilon}} \int_{Q''} \left( \fint_{B_\epsilon(x)} |u(x) - u(y)|\, dy \right) dx \leq C_2 (3\epsilon)^{N-1} \sum_{Q' \in \mathcal{Q}_{3\epsilon}} \fint_{Q'} \fint_{Q'} |u(x) - u(y)|\, dy\, dx \leq C_2 \kappa'_{3\epsilon}(u).$$

Proceeding by suitable translations of the family $\mathcal{Q}_{3\epsilon}$, we can make sure that each of the subcubes of $Q'$ is the inner cube of some cube in a translated family. Summing over the $3^N$ translated families of cubes, by additivity of the integral we then get

$$\frac{1}{\epsilon}\int_{\mathbb{R}^N}\left(\fint_{B_\epsilon(x)}|u(x)-u(y)|\,\mathrm{d}y\right)\mathrm{d}x \leq 3^N C_2 \kappa'_{3\epsilon}(u). \qquad \Box$$

For the same choice of mollifiers (5.1) in the functional of Theorem 1.1, one finds

$$(5.3)\quad \int_{\mathbb{R}^N}\left|\int_{\mathbb{R}^N}\frac{u(x)-u(y)}{|x-y|}\rho_\epsilon(x-y)\,\mathrm{d}y\right|\mathrm{d}x = \frac{C_N\omega_N}{\epsilon}\int_{\mathbb{R}^N}\left|u(x)-\fint_{B_\epsilon(x)}u(y)\,\mathrm{d}y\right|\mathrm{d}x.$$

Now, according to our results the energy (5.3) *cannot be comparable* to $\kappa_\epsilon$, $\kappa'_\epsilon$ or (5.2) on all scales for arbitrary $BV$ functions. Indeed, Theorem 1.1 implies that (5.3) tends to zero for $u \in W^{1,1}(\mathbb{R}^N)$ while $\kappa_\epsilon$, $\kappa'_\epsilon$ and (5.2) tend to a constant times the total variation of the weak derivative $\nabla u$ as $\epsilon$ tends to zero. They are nevertheless all *comparable* for $BV$ functions of the form $u = \chi_A$ in view of the identity

$$\int_{\mathbb{R}^N}\left|\int_{\mathbb{R}^N}\frac{\chi_A(x)-\chi_A(y)}{|x-y|}\rho_\epsilon(x-y)\,\mathrm{d}y\right|\mathrm{d}x = \int_{\mathbb{R}^N}\int_{\mathbb{R}^N}\frac{|\chi_A(x)-\chi_A(y)|}{|x-y|}\rho_\epsilon(x-y)\,\mathrm{d}y\,\mathrm{d}x.$$

ACKNOWLEDGEMENTS

The authors would like to thank the referee for his or her detailed reading and comments that have improved the quality of the presentation of the paper. This work was initiated while the first author (ACP) was visiting the Department of Applied Mathematics of the National Chiao Tung University (Research Grant 104W986). He warmly thanks the Department for the invitation and hospitality; he has also been supported by the Fonds de la Recherche scientifique–FNRS under research grant J.0026.15. The second author (DS) is supported by the Taiwan Ministry of Science and Technology under research grants 103-2115-M-009-016-MY2 and 105-2115-M-009-004-MY2.

AUGUSTO C. PONCE
UNIVERSITÉ CATHOLIQUE DE LOUVAIN
INSTITUT DE RECHERCHE EN MATHÉMATIQUE ET PHYSIQUE
CHEMIN DU CYCLOTRON 2, L7.01.02
1348 LOUVAIN-LA-NEUVE
BELGIUM
*E-mail address*: `Augusto.Ponce@uclouvain.be`

DANIEL SPECTOR
NATIONAL CHIAO TUNG UNIVERSITY
DEPARTMENT OF APPLIED MATHEMATICS
HSINCHU, TAIWAN
*E-mail address*: `dspector@math.nctu.edu.tw`